\documentclass[a4paper,12pt]{article}\usepackage[]{graphicx}\usepackage[]{color}
\makeatletter
\def\maxwidth{ %
  \ifdim\Gin@nat@width>\linewidth
    \linewidth
  \else
    \Gin@nat@width
  \fi
}
\makeatother

\definecolor{fgcolor}{rgb}{0.345, 0.345, 0.345}

\usepackage{framed}
\RequirePackage[numbers]{natbib}
\makeatletter
 {\par\unskip\endMakeFramed%
 \at@end@of@kframe}
\makeatother

\definecolor{shadecolor}{rgb}{.97, .97, .97}
\definecolor{messagecolor}{rgb}{0, 0, 0}
\definecolor{warningcolor}{rgb}{1, 0, 1}
\definecolor{errorcolor}{rgb}{1, 0, 0}

\usepackage{alltt}
\usepackage{amsfonts, amsmath, amsthm, amssymb}
\usepackage[english]{babel}
\usepackage{booktabs}
\usepackage{bm}
\usepackage{eucal}
\usepackage{enumerate}
\usepackage{float}
\usepackage{dsfont}
\usepackage{epsfig}
\usepackage{fontenc}
\usepackage[hmargin = 1.6cm, vmargin = 2.3cm]{geometry}
\usepackage{lscape}
\usepackage{mathrsfs}
\usepackage{natbib}
\usepackage{rotating}
\usepackage{verbatim}
\usepackage[table]{xcolor}

\newcommand{\dif}{\ensuremath{\mathrm{d}}}

\newcommand{\T}{\ensuremath{\mathrm{\scriptscriptstyle T}}}
\usepackage{bm}

\usepackage{amsmath}

\DeclareMathOperator{\E}{E}

\usepackage{mathrsfs}

\newcommand{\vs}{\vspace{0.2cm}}

\newtheorem{theorem}{Theorem}

\newtheorem{example}{Example}

\newtheorem{corollary}[theorem]{Corollary}

 \let\oldthebibliography=\thebibliography
 \let\oldendthebibliography=\endthebibliography
 \renewenvironment{thebibliography}[1]{%
   \footnotesize
   \oldthebibliography{#1}%
   \setlength{\parskip}{0mm}%
   \setlength{\itemsep}{0mm}%
 }{\oldendthebibliography}

\date{}

\renewcommand{\thefootnote}{$\dagger$}
\IfFileExists{upquote.sty}{\usepackage{upquote}}{}
\begin{document} 
\title{\vspace{-2cm}On the Kolmogorov Superposition Theorem \\ and Regular Means}
\author{Miguel de Carvalho}
\date{}
\maketitle 

\begin{abstract}\footnotesize
  While Kolmogorov’s probability axioms are widely recognized, it is less well known that in an often-overlooked 1930 note, Kolmogorov proposed an axiomatic framework for a unifying concept of the mean---referred to as \emph{regular means}. This framework yields a well-defined functional form encompassing the arithmetic, geometric, and harmonic means, among others.

In this article, we uncover an elegant connection between two key results of Kolmogorov by showing that the class of regular means can be derived directly from the Kolmogorov superposition theorem. This connection is conceptually appealing and illustrates that the superposition theorem deserves wider recognition in Statistics---not only because of its link to regular means as shown here, but also due to its influence on the development of neural models and its potential connections with other statistical frameworks. In addition, we establish a stability property of regular means, showing that they vary smoothly under small perturbations of the generator. Finally, we provide insights into a recent universal central limit theorem that applies to the broad class of regular means. 
\end{abstract}
\let\thefootnote\relax\footnotetext{M.~de Carvalho; School of Mathematics, University of Edinburgh; EH9 3FD UK.} \vspace{-0.5cm}

\section{Introduction}\label{intro}
\subsection{Background and Motivation}
While Kolmogorov's probability axioms are well known, his often-overlooked 1930 note \citep{kolmogorov1930} introduced an axiomatic framework for a unifying concept of the mean. Using it, Kolmogorov proved that the only functional form compatible with these axioms includes, as particular cases, the arithmetic, geometric, and harmonic means---among others. 

In this paper we build on a powerful but seemingly unrelated result established by Kolmogorov over 20 years later---commonly referred to as the Kolmogorov superposition theorem \cite[e.g.,][]{braun2009, lorentz1996}. When first proved, this theorem had major implications, as it refuted a conjecture implicit in Hilbert's 13th problem presented in 1900 at the International Congress of Mathematicians. See \cite{morris2021} for an historical account.

As the name suggests, the theorem concerns superpositions. In mathematical terms, a superposition is a function of functions; for example,
\begin{equation}\label{kol}
  f(x_1, x_2, x_3) = g\big(a(\alpha(x_1), \beta(x_2, x_3)),\, b(x_2, x_3)\big),
\end{equation}
is a superposition of univariate and bivariate functions. The superposition theorem states that any \textit{multivariate} continuous function can be expressed as sums and superpositions of \textit{univariate} continuous functions. In loose terms, any multivariate continuous function can therefore be represented exactly using only a finite number of univariate continuous functions. 

The influence of this theorem persists in modern times, serving as the foundation for a wide range of neural models inspired by its principles \citep[e.g.,][]{decarvalho2025, fakhoury2022, lin1993, liu2024, montanelli2020, schmidt2021, sprecher2002}. These neural models interpret the two distinct sets of univariate continuous functions in the superposition as two different layers: the so-called inner functions correspond to the mapping in the hidden layer, while the outer functions correspond to the mapping in the output layer. The superposition theorem is the theoretical foundation of such models, in the same way that multi-layer perceptrons are grounded on the universal approximation theorem.  


As shown in this paper, these two seemingly unrelated results by Kolmogorov---i.e., functional form of regular means and the  superposition theorem---have an important link that has been underappreciated. Specifically, the functional form of regular means can be derived as a corollary of the Kolmogorov superposition theorem. This connection is conceptually appealing and illustrates that the superposition theorem deserves wider recognition in Statistics---not only because of its link to regular means, as shown here, but also due to its influence on the development of neural models, as noted earlier, and its potential connections with other statistical frameworks, as will be hinted at later. In addition, we establish a stability property of regular means, showing that they vary smoothly under small perturbations of the generator. We relate this result to the axiomatics of regular means, and to approximations of the geometric mean---which play a fundamental role in modern portfolio theory \citep[][Chapter~3]{markowitz2013}.

Finally, we provide insights into a recent universal central limit theorem that applies to the broad class of regular means. This result, originally due to \citet{decarvalho2016a}, establishes a central limit theorem valid for all types of means, including the arithmetic, geometric, and harmonic means. Despite its generality, little was previously known about its finite-sample behavior, both theoretically and empirically---a gap that this work aims to address. We offer insights into the factors affecting the speed of convergence to the normal approximation via an Edgeworth expansion and illustrate these findings through numerical experiments.
    

\subsection{Preparations on Regular Means}\label{regme}
Before proceeding, we first lay the groundwork for regular means. Let $\mathbf{x} = (x_1, \dots, x_n)^{\T}$. 
A \textit{regular mean} is a mapping $M:\mathbb{R}^n \to \mathbb{R}$ satisfying the following axioms:
    \begin{enumerate}
    \item[-\,A1)] $M(\mathbf{x})$ is continuous and increasing in each variable. \vspace{-0.2cm}
    \item[-\,A2)] $M(\mathbf{x})$ is a symmetric function. \vspace{-0.2cm}
    \item[-\,A3)] $M(x\mathbf{1}_n) = x$, i.e., the mean of $n$ identical values equals their common value.\vspace{-0.2cm}
    \item[-\,A4)] The mean of $\mathbf{x}$ is unchanged if any component is replaced by its mean, 
 $m= M(x_1, \ldots, x_{n_0})$, i.e. $M(\mathbf{x}) = M(\underbrace{m, \ldots, m}_{n_0}, x_{n_{0}+1}, \ldots, x_{n})$.
    \end{enumerate}
\noindent \cite{kolmogorov1930} proved that if conditions A1 to A4 are satisfied, the function $M(\mathbf{x})$ takes the form 
    \begin{equation}
      M(\mathbf{x}) \equiv M_g(\mathbf{x}) = g^{-1} \bigg(\frac{1}{n}\sum_{i=1}^n g(x_i)\bigg), 
      \label{rmean}
    \end{equation}
    where $g: \Omega \subseteq\mathbb{R} \to \mathbb{R}$ is a continuous monotone function and $g^{-1}$ denotes its inverse. For reasons that will become clear later, we refer to $g$ as the generator (or activation) function and to $\Omega$ as the input range. For example, if $g(x) = x$ with $\Omega = \mathbb{R}$, we obtain the arithmetic mean,
$$M_x(\mathbf{x}) \equiv \frac{1}{n} \sum_{i = 1}^n x_i.$$  
Also, $g(x) = \log x$ and $g(x) = 1/x$, with $\Omega = (0, \infty)$, correspond respectively to the geometric and harmonic means:
\begin{equation*}
  \begin{split}
    M_{\log \, x}(\mathbf{x}) = \left(\prod_{i = 1}^n x_i\right)^{1/n}, \qquad 
    M_{1/x}(\mathbf{x}) = \frac{n}{\sum_{i = 1}^n x_i^{-1}}.
    \end{split}\vspace{-0.4cm}
\end{equation*}
Equation~\eqref{rmean} also includes as particular cases other examples, such as the so-called power mean: For example, if $g(x) = x^p$, for $p \in (0, \infty)$, we obtain $M_{x^p}(\mathbf{x}) = \{1/n \sum_{i=1}^n x_i^p\}^{1/p}$.

    Kolmogorov's proof of the result in \eqref{rmean} is far from straightforward,
    requiring the use of an auxiliary function and a distinct approach for
    rational and irrational numbers. See 
    \citet[][Chapter~5]{leinster2021} for an alternative approach.

\subsection{Outline}
Section~\ref{super} introduces the superposition theorem and derives two corollaries, including one that establishes \eqref{rmean}. In light of this connection, Section~\ref{generators} discusses properties of the generators, while Section~\ref{continuity} examines the continuity of regular means. Section~\ref{UnivCLT} presents the universal central limit theorem, provides rate insights via an Edgeworth expansion, and offers numerical illustrations of the result.

\section{On Regular Means as Superpositions}\label{superpositions}
\subsection{Means are Superpositions}\label{super}
Below I present a modern version of the Kolmogorov superposition theorem, equivalent to that in  \cite{lorentz1996}. 
\begin{theorem}\label{KLAT}
  Let $f:A \rightarrow \mathbb{R}$ be a multivariate continuous function over a compact set $A \subset \mathbb{R}^n$. Then, it has the representation
  \begin{equation}\label{KLAT2}
    f(\mathbf{x}) = \sum_{j = 1}^{2 n + 1}l\left(\sum_{i = 1}^n \lambda_ig_j(x_i) \right), 
  \end{equation}
  with continuous one–dimensional outer and inner functions $l$ and $g_{j}$, where $g_{j}$ is increasing, and  $\sum_{i = 1}^n \lambda_i \leq 1$ with $\lambda_i \geq 0$, for all $i$.
\end{theorem}

\noindent In simple terms, any continuous multivariate function can be expressed as a function of univariate functions. In mathematical parlance, such functions of functions are referred to as superpositions. Thus, Theorem~\ref{KLAT} can be summarized as stating that any continuous multivariate function can be represented using sums and superpositions of univariate functions. It is typical to set $A = [0, 1]^n$, but it is well-known that the superposition theorem holds more generally for compact metric spaces \citep[][]{ostrand1965}. 

The functional form of the regular mean \eqref{rmean} resembles that of the superposition theorem in \eqref{KLAT2}. In fact, by setting in \eqref{KLAT2},  
\begin{equation*}
  \lambda_i = \frac{1}{n}, \quad g_j(x) = g(x), \quad l(x) = \frac{1}{2n + 1} g^{-1}(x),
\end{equation*}
we recover the regular mean in \eqref{rmean}.

From all axioms that a regular mean needs to obey, only continuity is required by Theorem~\ref{KLAT}. Adding symmetry (i.e., A2) to the list of assumptions of Theorem~\ref{KLAT} leads to the following result whose proof is given in the Appendix; a comparable result is asserted in \cite{zaheer2017}.

\begin{corollary}\label{KLATB}
  Let $f:A \rightarrow \mathbb{R}$ be a symmetric multivariate continuous function over a compact set $A \subset \mathbb{R}$. Then it has the representation
  \begin{equation}\label{KLATB2}
    f(\mathbf{x}) = l\left(\lambda \sum_{i = 1}^n g(x_i) \right), 
  \end{equation}
  with continuous one–dimensional outer and inner functions $l$ and $g$, where $g$ is increasing and $\lambda n \leq 1$ with $\lambda \geq 0$.
\end{corollary}

Equation~\eqref{KLATB2} is already a step closer to the functional form of the regular mean in \eqref{rmean}. Next, we show that A1, A2, A3, and A4 finally provide us with all we need to reach \eqref{rmean}. 


\begin{corollary}\label{KLATC}
  Let $f:A \rightarrow \mathbb{R}$ be a multivariate function obeying A1 to A4. Then it has the representation
  \begin{equation}\label{KLATC2}
    f(\mathbf{x}) = g^{-1}\left(\frac{1}{n}\sum_{i = 1}^n g(x_i) \right), 
  \end{equation}
  with continuous one–dimensional and increasing function $g$.
\end{corollary}
\noindent While this result is well known \citep[e.g.,][]{decarvalho2016a}, its derivation from the superposition theorem (Theorem~\ref{KLAT}) is entirely novel and uncovers an elegant connection between two fundamental results of Kolmogorov. In words, Corollary~\ref{KLATC} provides a superposition theorem for regular means, showing that the only functional form satisfying A1--A4 is \eqref{KLATC2}; that is, the only admissible functions obeying A1--A4 are regular means.

\subsection{Regular Means as Simple Neural Models?}\label{generators}

Since it follows from Section~\ref{super} that regular means share the same structure as~\eqref{KLAT2}, we next draw a parallel with the neural models based on that representation. Before proceeding, we emphasize that regular means should \textit{not} be regarded as actual neural models, although the analogy is instructive. Under this interpretation, $g$ and $g^{-1}$ play the roles of activation functions,
with only continuous and monotone functions being admissible choices for $g$---thus excluding, for instance, nonlinearities such as GELU (Gaussian Error Linear Unit). Table~\ref{activation_means} lists common means and their corresponding activation functions. While some of these functions are uncommon in deep learning, they arise naturally in the context of Kolmogorov’s regular means.

  \begin{table*}[h]
\caption{Common means, their corresponding activation functions, and Kolmogorov expected values as defined in \eqref{kev}.}
\label{activation_means}\footnotesize \centering \vs 
\begin{tabular}{@{}lllll@{}}
\hline
  \textbf{Activation} & \textbf{Input} & \textbf{Resulting} & \textbf{Expression} & \textbf{Kolmogorov} \\
  \textbf{function} ($g$) & \textbf{range} ($\Omega$) & \textbf{mean} & ($M$) & \textbf{expected value} \\ 
\hline
$x$ & $\mathbb{R}$            & Arithmetic & $1/n \sum_{i = 1}^n x_i$ & $\E(X)$\\
  $\log(x)$ & $x > 0$            & Geometric & $(\prod_{i = 1}^n x_i)^{1/n}$ & $\exp\{\E(\log X)\}$\\  
$1/x$ & $x > 0$               & Harmonic & $n/\sum_{i = 1}^n x_i^{-1}$ & $1 / \E(1/X)$ \\
  $x^p$ & $x > 0$               & Power & $\{1/n \sum_{i=1}^n x_i^p\}^{1/p}$ & $\{\E(X^p)\}^{1/p}$ \\
  $\exp(x)$ & $\mathbb{R}$ & Exponential & $\log\{1/n \sum_{i = 1}^n \exp(x_i)\}$ & $\log[\E\{\exp(X)\}]$ \\
\hline
\end{tabular}
\end{table*}

\subsection{On the Continuity of Regular Means}\label{continuity}
The main theoretical insight of this paper, established in Section~\ref{super}, is that the functional form of all means in~\eqref{rmean} follows as a corollary of the Kolmogorov superposition theorem. Although secondary, the next result is of independent interest. Let $C(A)$ denote the space of continuous functions on a compact set $A \subset \mathbb{R}^n$, and let $C_{\mathrm{LI}}(B)$ denote the space of Lipschitz continuous, strictly monotone functions on a compact set $B \subset \Omega$. Both spaces are endowed with the supremum norm
\[
\|f\|_{\infty} = \sup_{x \in X} |f(x)|, \qquad X \in \{A,B\}.
\]
The following result holds.
\begin{theorem}\label{cont}
  Let
  $M: C_{\emph{LI}}(B) \to C(A)$ be the operator 
  defined as $M_g(\mathbf{x}) = g^{-1}(\sum_{i = 1}^ng(x_i))$ where $A$ is a compact subset of $\mathbb{R}^n$, and $B$ is a compact subset of $\Omega$. Then, $M_g$ is continuous in $g$.
\end{theorem}
Whereas continuity of the regular mean as a function of the data $\mathbf{x}$ is required by A1, Theorem~\ref{cont} establishes that $M_g$ is itself continuous as a function of the generator $g$. This property guarantees the stability of the regular mean, ensuring that the regular mean does not change abruptly under small perturbations of $g$. Without this continuity, two nearly identical activation functions $g$ and $h$ could yield markedly different means, making $M_g(\mathbf{x})$ unstable or unreliable. The conditions under which the result is derived are mild; for instance, all activation functions listed in Table~\ref{activation_means} are Lipschitz continuous on any compact set~$B \subset \Omega$.

Theorem~\ref{cont} indicates that one can move from one mean based on a generator to another defined by a nearby generator, with the two means remaining close. This stability principle is widely used in finance as shown in the following example.

\begin{example}[Portfolio theory]\normalfont
The geometric mean is the appropriate notion of average return for a portfolio with period returns $r_1,\dots,r_n$. Indeed, if $w_0$ is the initial wealth, then
\begin{equation*}
  \begin{split}
w_n = w_0 (1+r_1)\cdots(1+r_n)
     = w_0 \,\{M_{\log x}(1+r_1,\dots,1+r_n)\}^n,
  \end{split}
\end{equation*}
so the average return over $n$ periods is the geometric mean
\[
M_{\log x}(1+r_1,\dots,1+r_n)
  = \Bigg\{\prod_{t=1}^n (1+r_t)\Bigg\}^{1/n}.
\]
Since portfolio theory is typically formulated in a mean-variance framework, the geometric mean is often approximated by a first-order Taylor expansion around the sample mean $\bar r$, giving
\[
M_{\log x}(1+r_1,\dots,1+r_n)
  \approx \exp\{\bar r - (\bar r^2 + s^2)/2\},
\]
where $\bar r$ and $s^2$ denote the sample mean and variance of $r_1,\dots,r_n$. For this and related approximations and their applications see  \citet[][]{markowitz2014}.
\end{example}

\section{On a Universal Central Limit Theorem}\label{UnivCLT}
\subsection{Large Sample Behavior of Regular Means}\label{ls}
While the standard central limit theorem is about the limiting behavior of the arithmetic mean, Theorem~\ref{UCLT} is a general central limit theorem that applies to the entire class of regular means, including those reported in Table~\ref{activation_means}.
Throughout, let $\textbf{X} = (X_1, \ldots, X_n)$ and define the Kolmogorov expected value, as
\begin{equation}\label{kev}
  E_g (X) = g^{-1} (E_x\{g(X)\}), 
\end{equation}
where $E_x\{g(X)\} = \int g(x) \, \dif F(x)$, with $X \sim F$, and $g$ is a continuous monotone function
\begin{theorem}[Universal Central Limit Theorem]\label{UCLT}
Suppose that $X_1, 
\ldots, X_n \overset{\emph{iid}}{\sim} F_X$, and that $\emph{var}\{g(X)\} < \infty$, where $g(x)$ is a strictly monotone function with derivative $g'(E_g(X))$ at $x = E_g(X)$. Then, as $n \to \infty$, 
\begin{equation}\label{avar}
  \sqrt{n}\{M_g(\mathbf{X}) - E_g(X)\}
  \overset{\emph{d}}{\rightarrow} \emph{N}\bigg(0,
  \frac{\emph{var}\{g(X)\}}{\{g'(E_g(X))\}^2}\bigg), 
\end{equation}
with $M_g(\mathbf{X})$ as in \eqref{rmean}.
\end{theorem}
\noindent The universal central limit theorem is due to de Carvalho \cite{decarvalho2016a}. While the result is well established, its convergence properties and numerical performance remain only partially understood. In the remainder of this section, we aim to shed further light on these issues.

To gain insight into the factors influencing the speed of convergence to the normal distribution, we employ Edgeworth expansions. Let $F_n$ denote the distribution function of
\begin{equation*}
\sqrt{n}\bigg(\frac{M_g(\mathbf{X}) - E_g(X)}
{\sqrt{\text{var}\{g(X)\}}/g'(E_g(X))}\bigg),
\end{equation*}
and let $G_n$ denote the distribution function of
\begin{equation*}
\sqrt{n}\bigg(\frac{n^{-1}\sum_{i=1}^n g(X_i) - E\{g(X)\}}
{\sqrt{\text{var}\{g(X)\}}}\bigg).
\end{equation*}

\noindent A standard delta-method argument \citep[][Section~2.7]{hall2013} 
shows that $F_n(x) = G_n(x) + o(1)$ as $n \to \infty$. 
Since the Edgeworth expansion for $G_n(x)$ is well known \citep[][Section~2.6]{hall2013}, it follows that
\begin{equation*}
  F_n(x) = \Phi(x) 
  - \phi(x)\bigg(\frac{\gamma_gp_1(x) }{6\sqrt{n}}
  + \frac{\kappa_gp_2(x) }{24n}
  + \frac{\kappa_g^2p_3(x)}{72n}\bigg)
  + o(n^{-1}),
\end{equation*}
where $\Phi(x)$ and $\phi(x)$ are respectively the distribution and density functions of the standard normal law, $\gamma_g$ and $\kappa_g$ denote respectively the skewness and excess kurtosis of $g(X)$, and $p_1$, $p_2$, and $p_3$ are Hermite polynomials (i.e., $p_1(x) = x^2 - 1$, $p_2(x) = x^3 - 3 x$, $p_3(x) = x^5 - 10 x^3 + 15 x$). The above Edgeworth expansion shows that heavier tails or greater skewness of $g(X)$ slow the convergence to the limiting normal distribution. As expected, when $g(x)=x$, the expansion reduces to the standard Edgeworth expansion for the arithmetic mean.



\begin{figure}[H]
   \includegraphics[scale = .055]{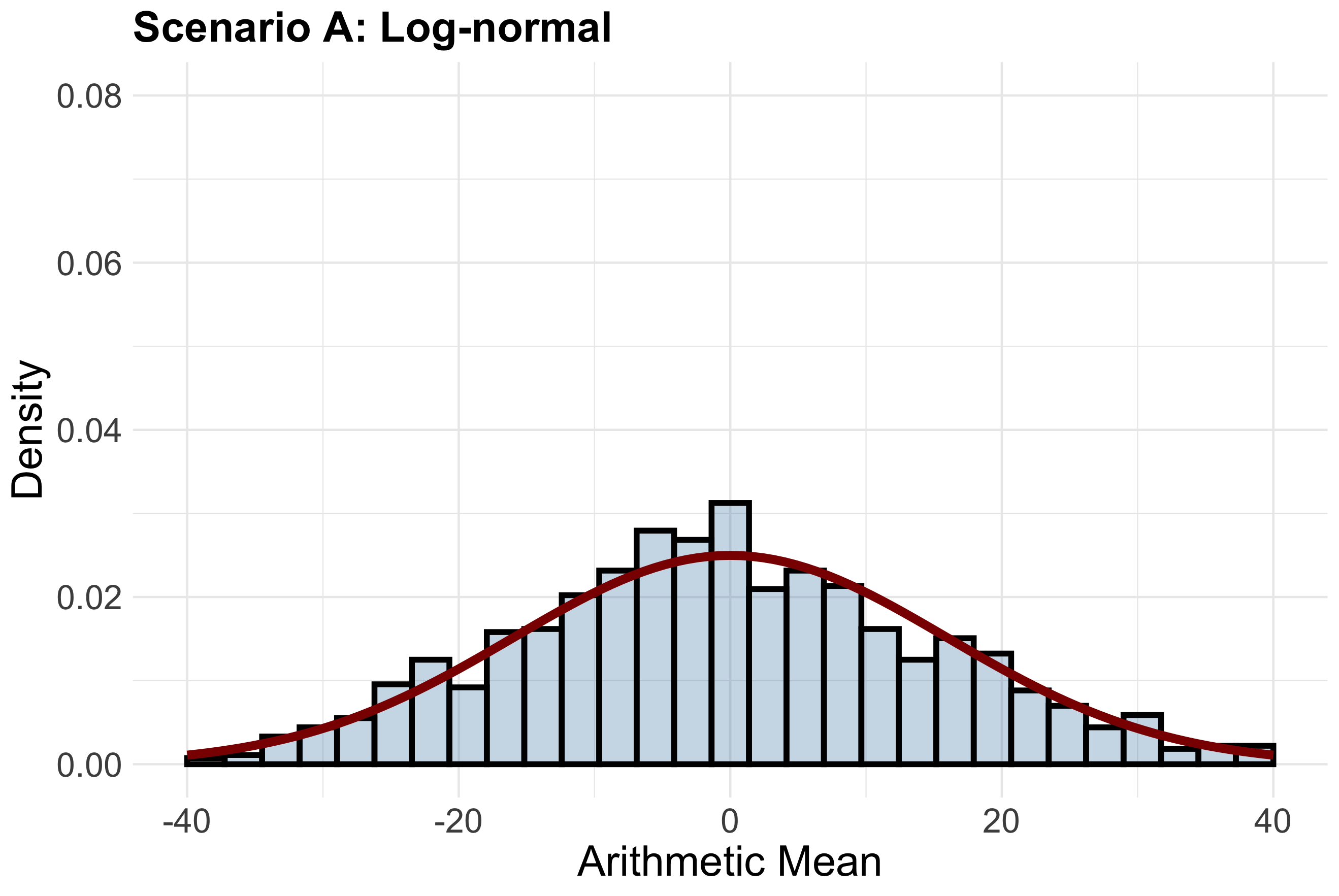}   \includegraphics[scale = .055]{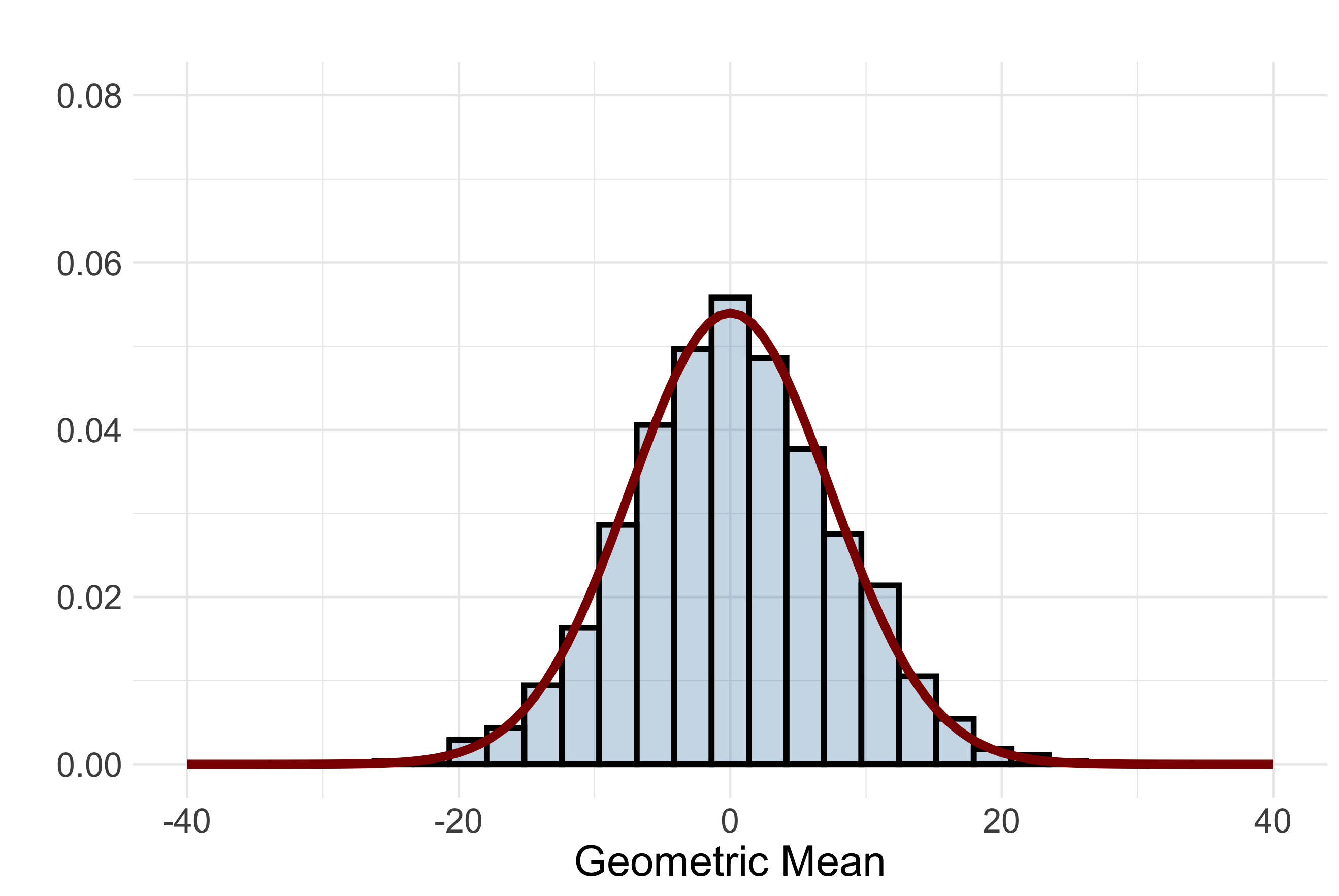} \includegraphics[scale = .055]{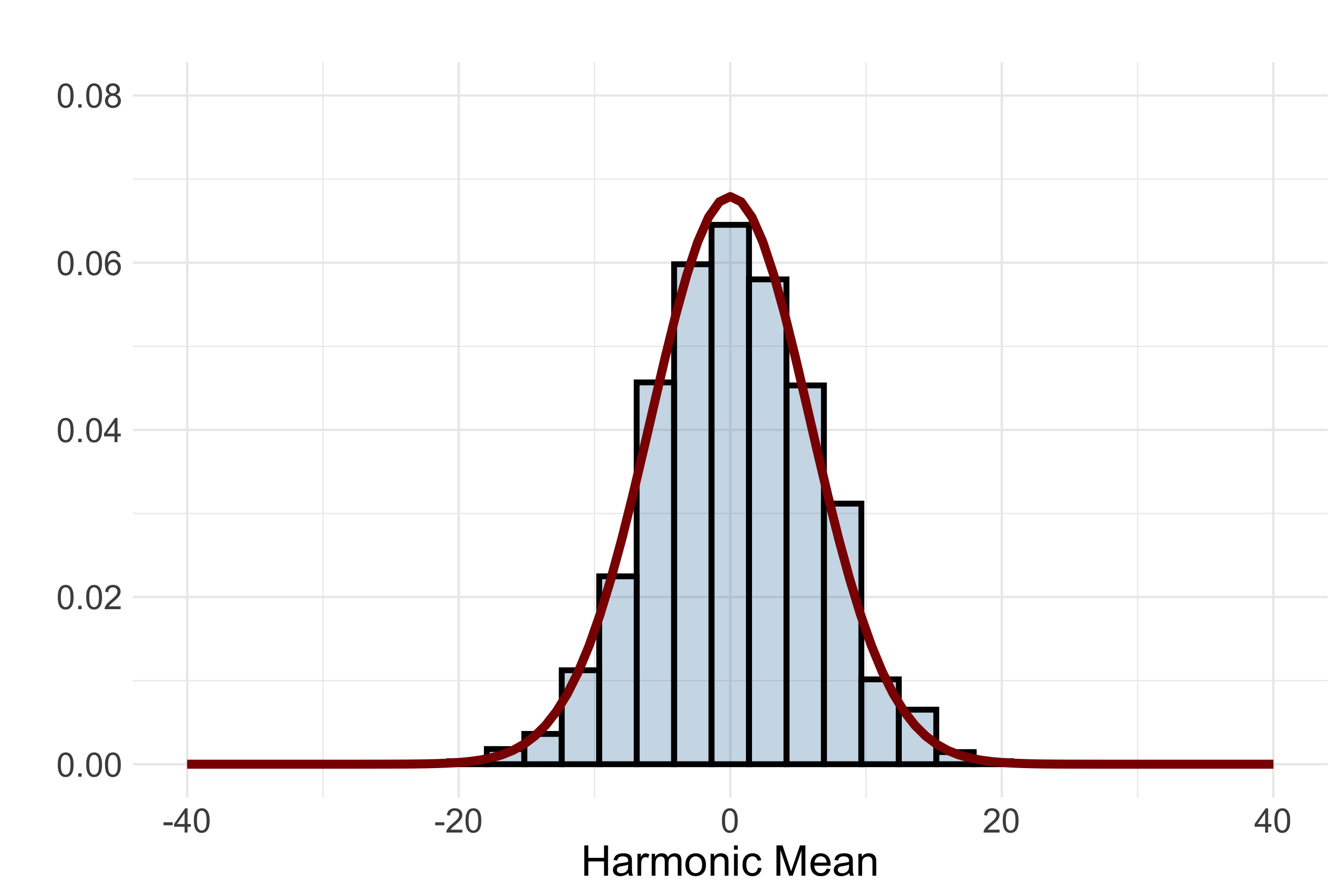}
  \includegraphics[scale = .055]{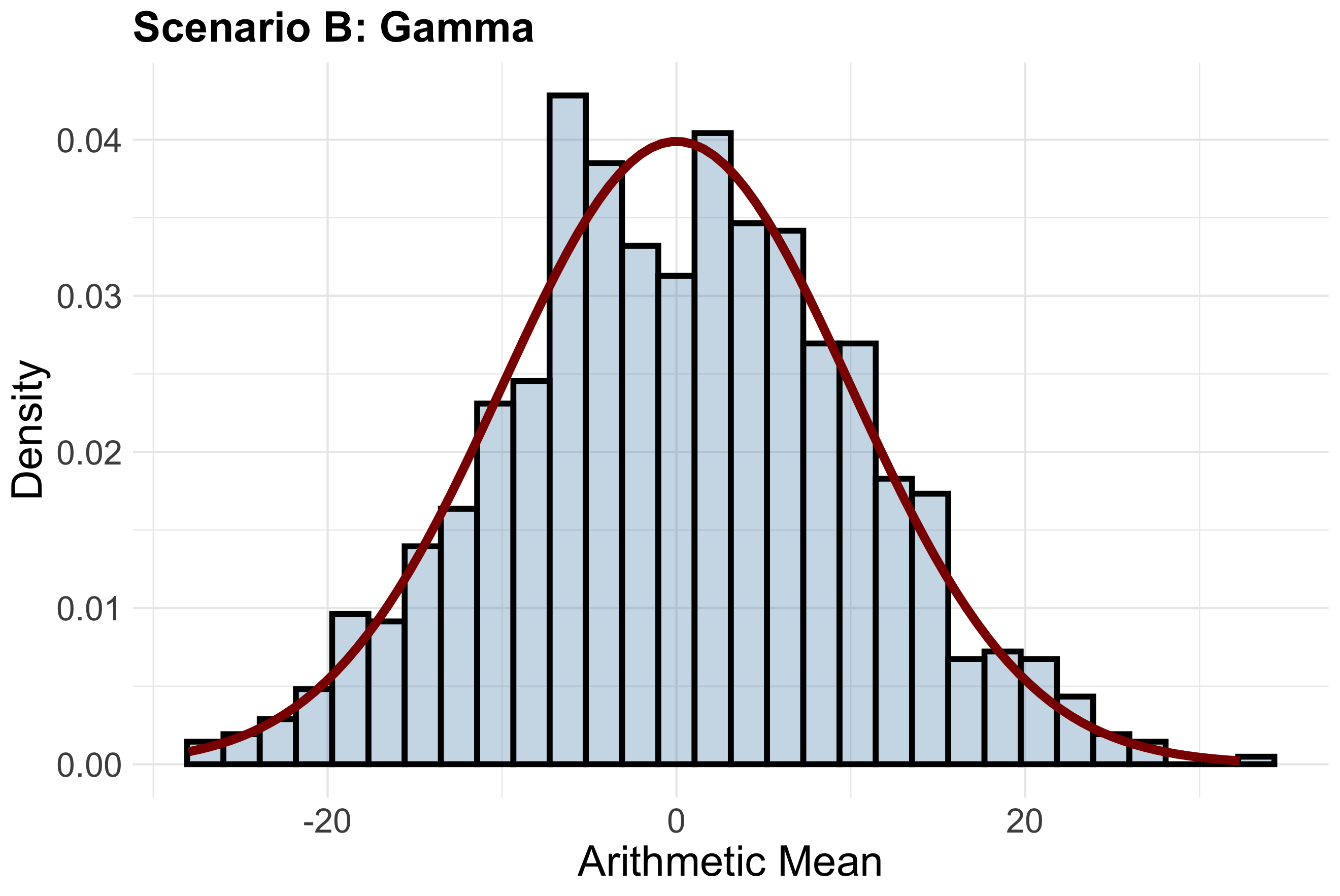}   \includegraphics[scale = .055]{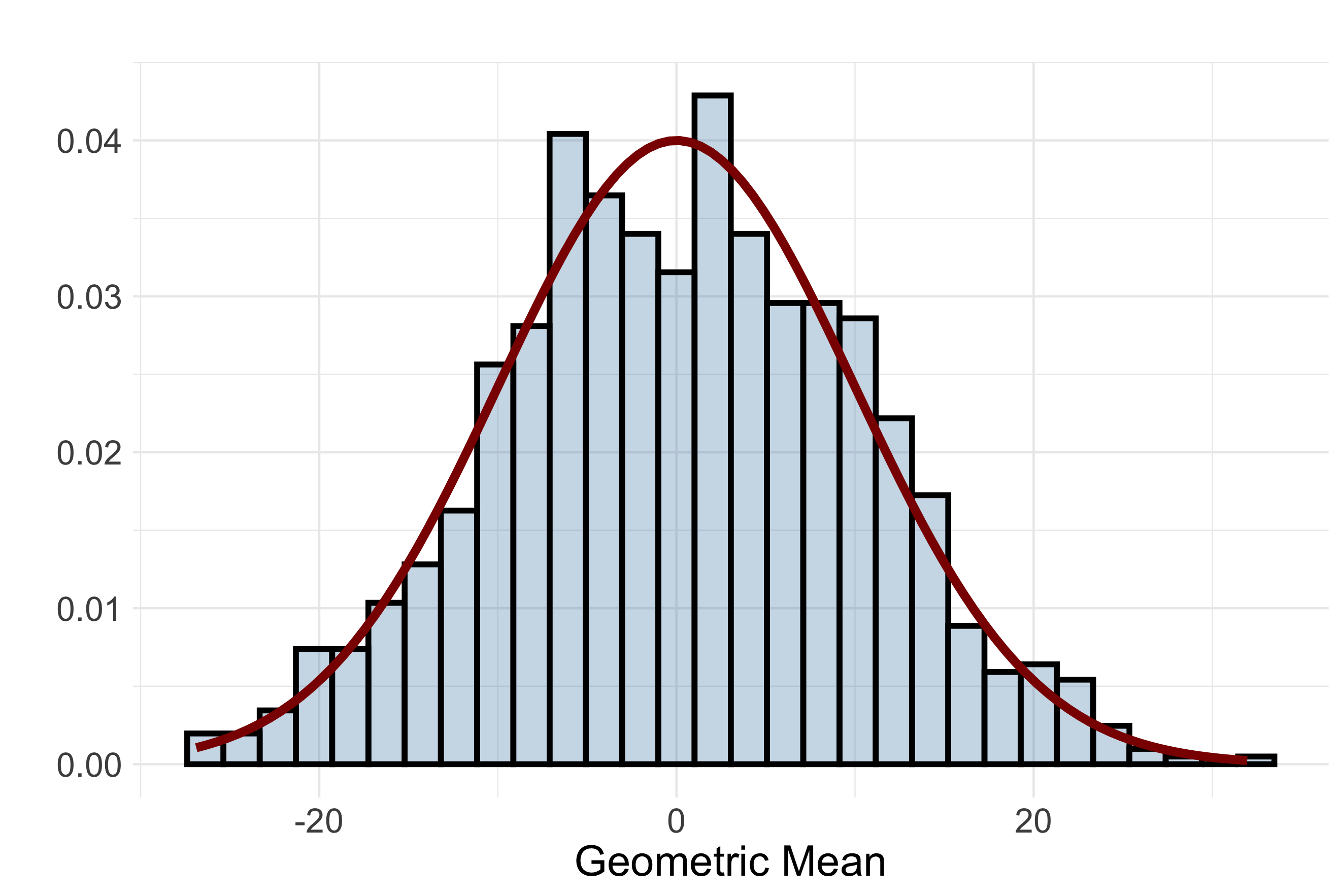} \includegraphics[scale = .055]{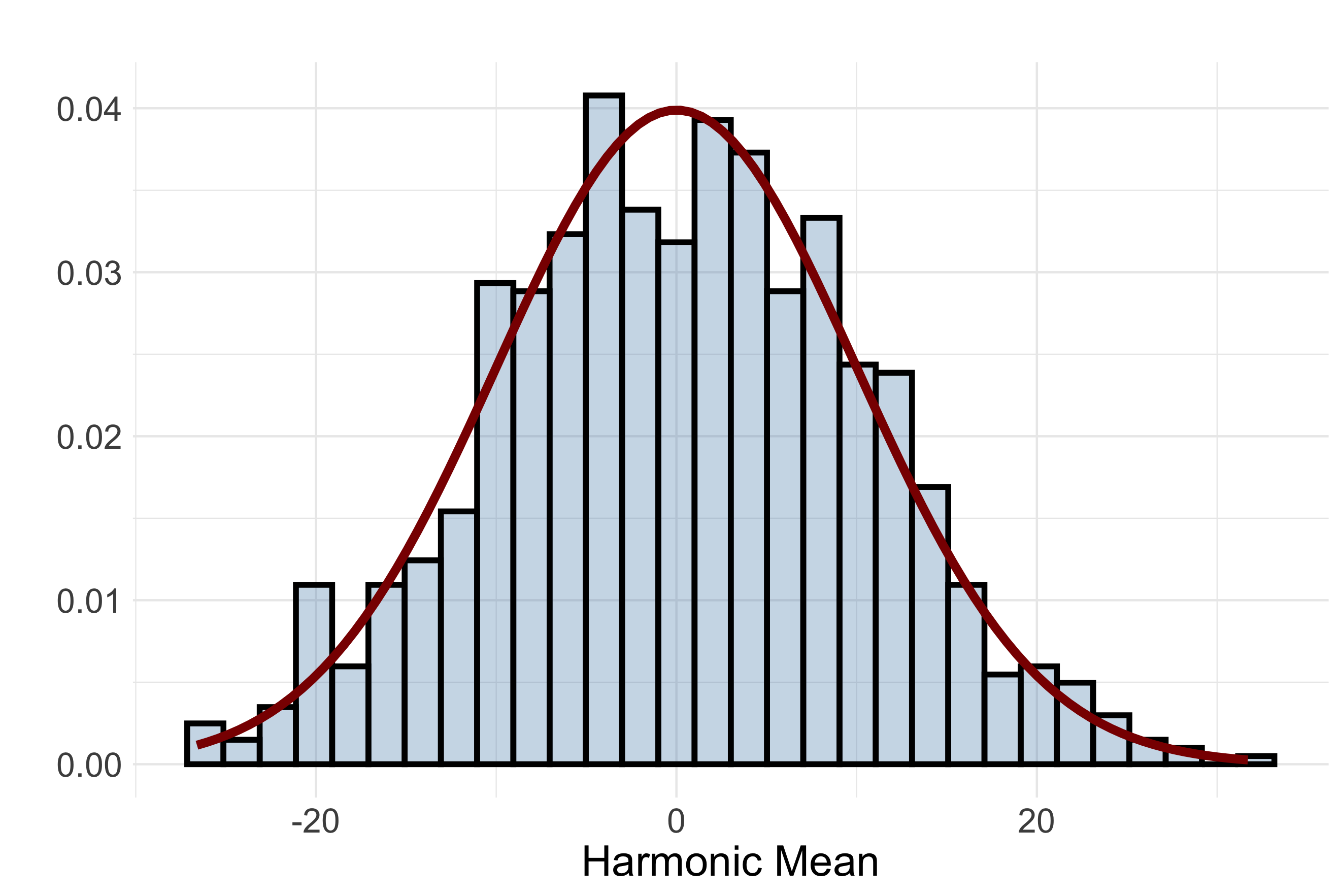}
  \includegraphics[scale = .055]{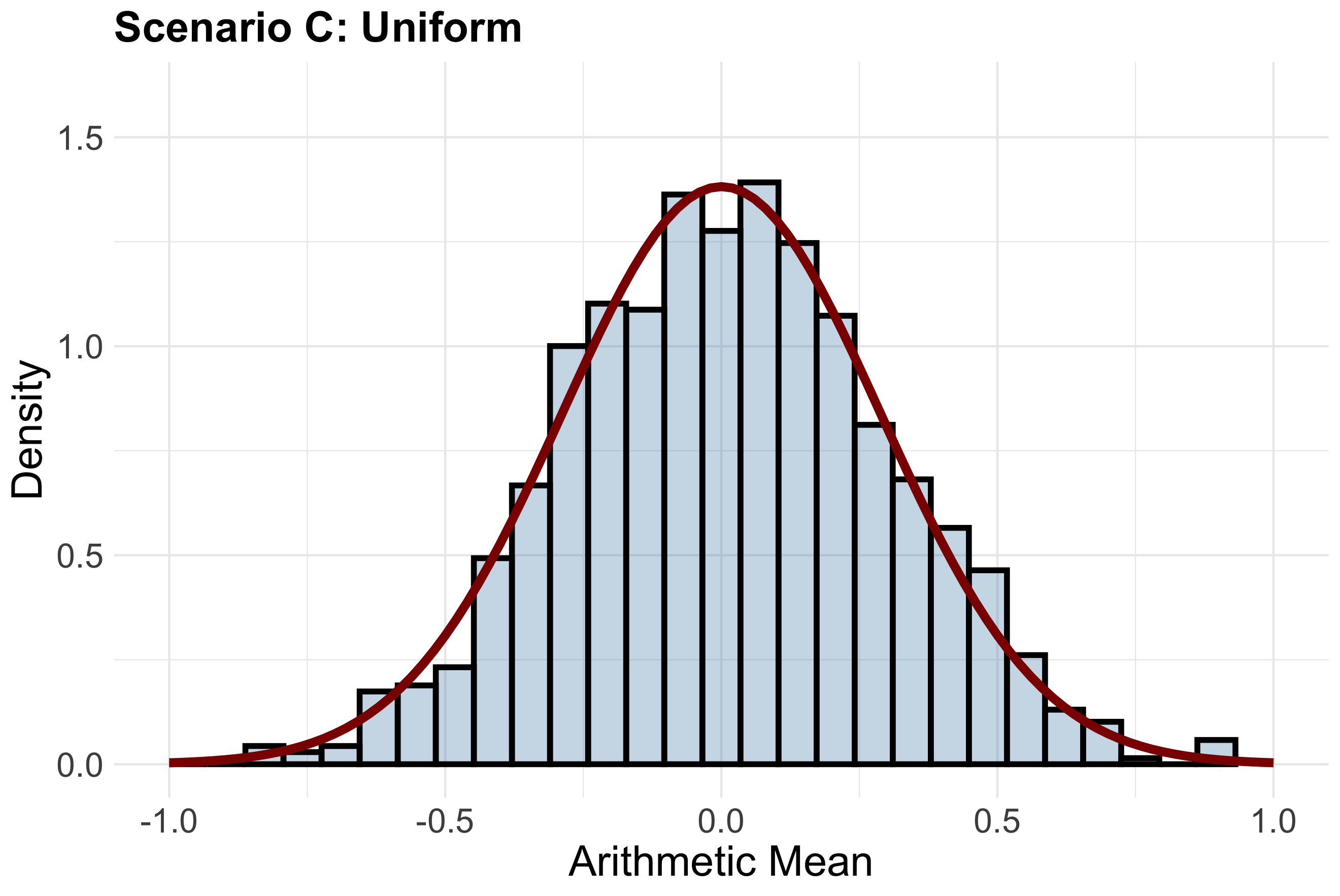} \includegraphics[scale = .055]{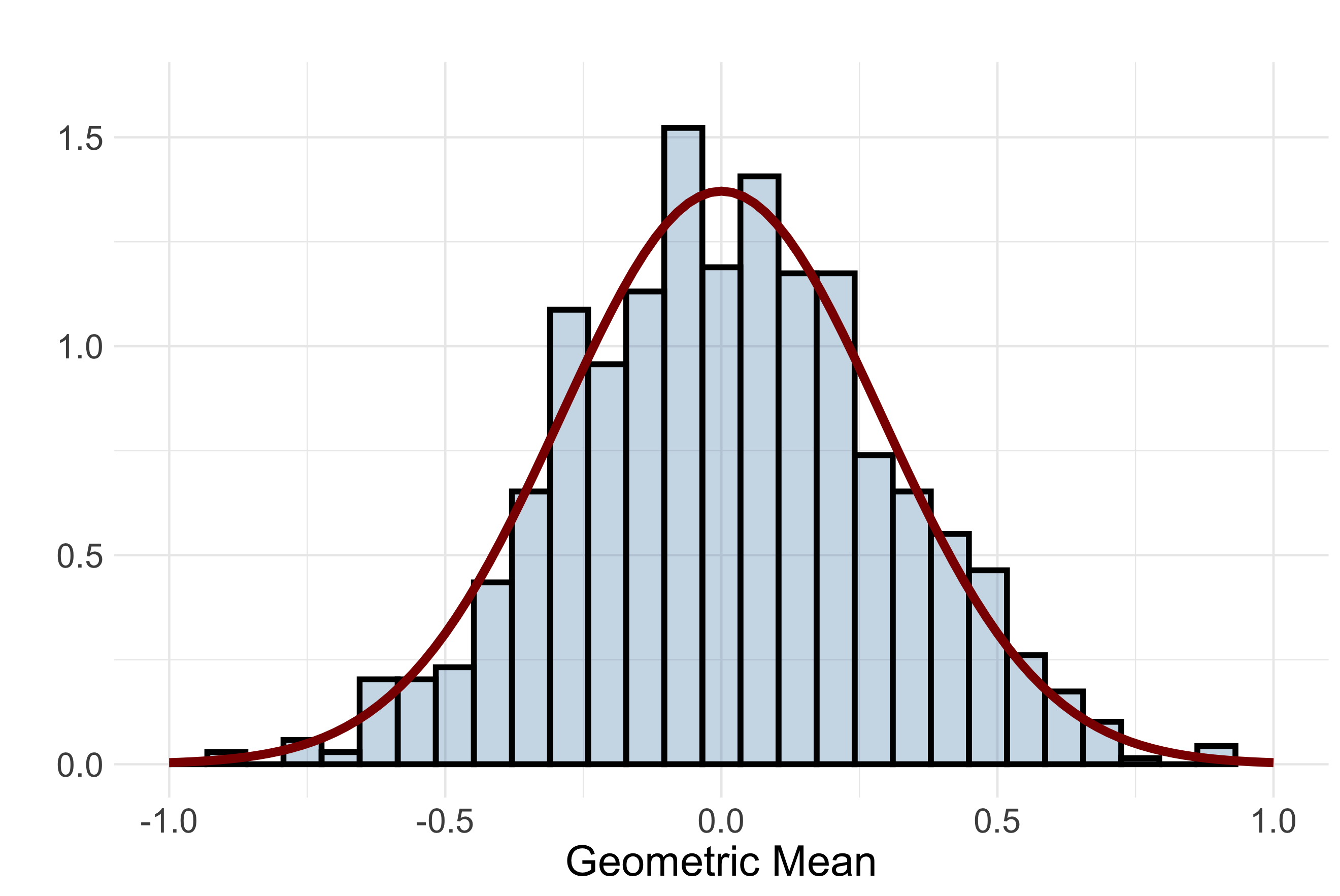} \includegraphics[scale = .055]{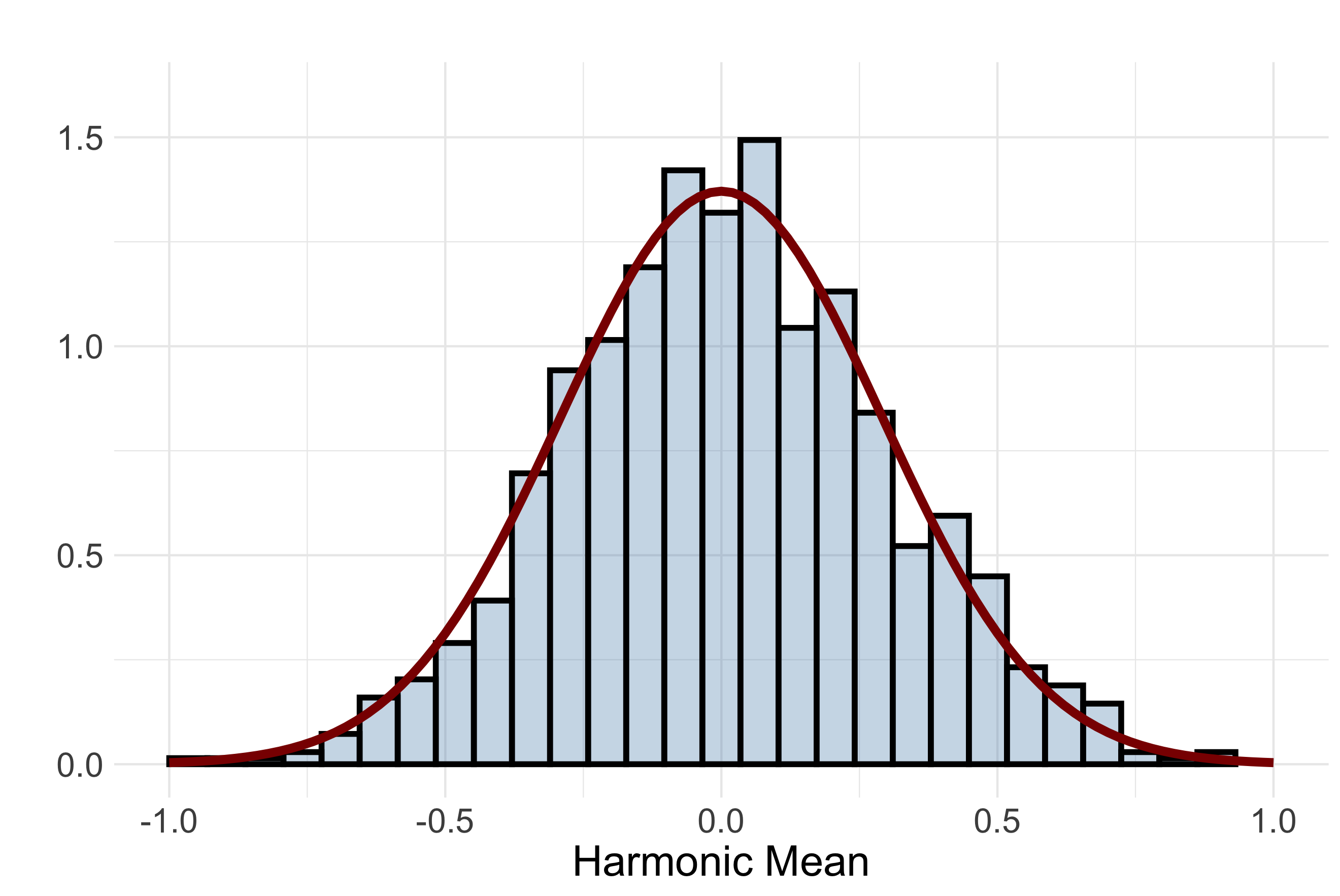}
    \includegraphics[scale = .055]{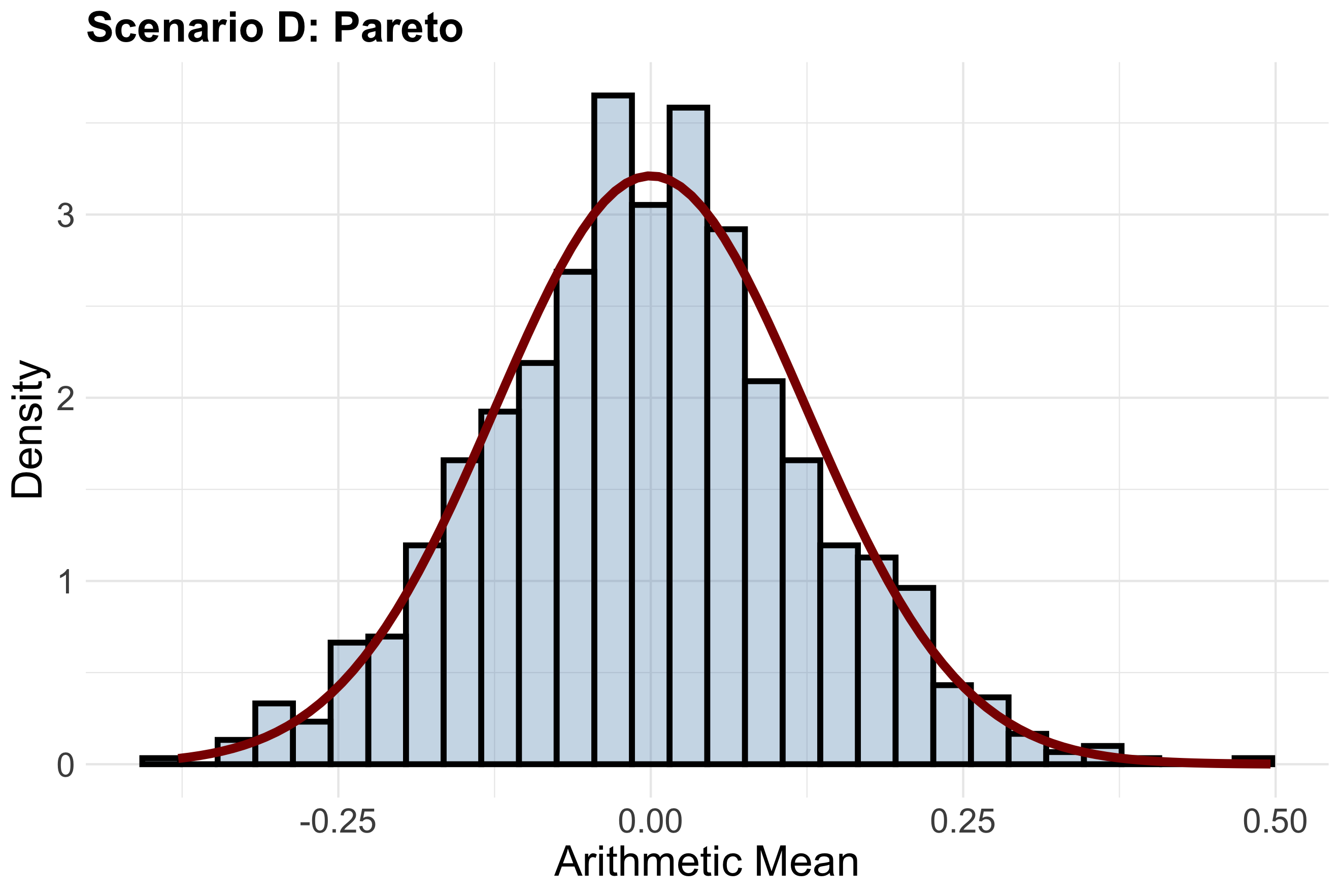} \includegraphics[scale = .055]{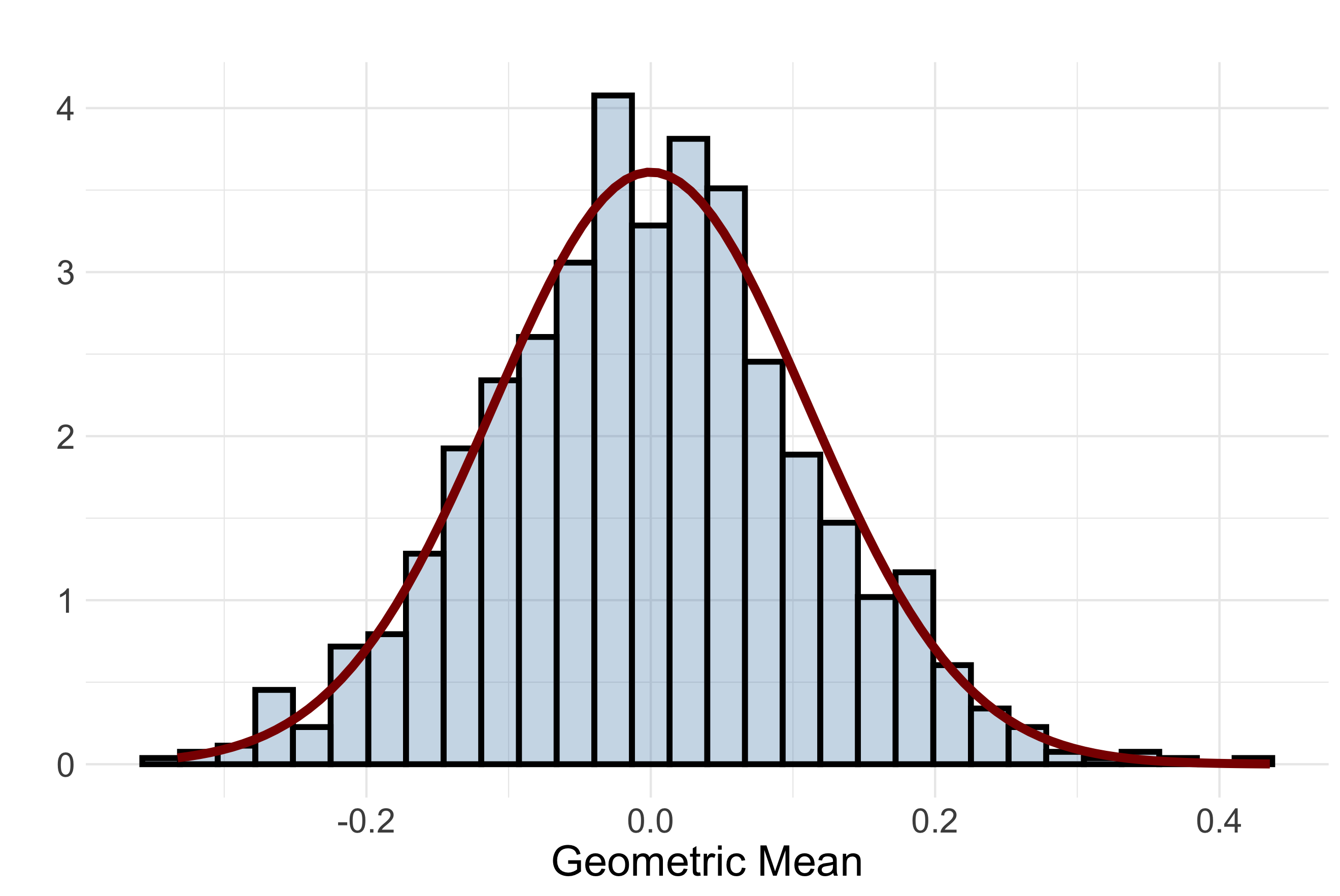} \includegraphics[scale = .055]{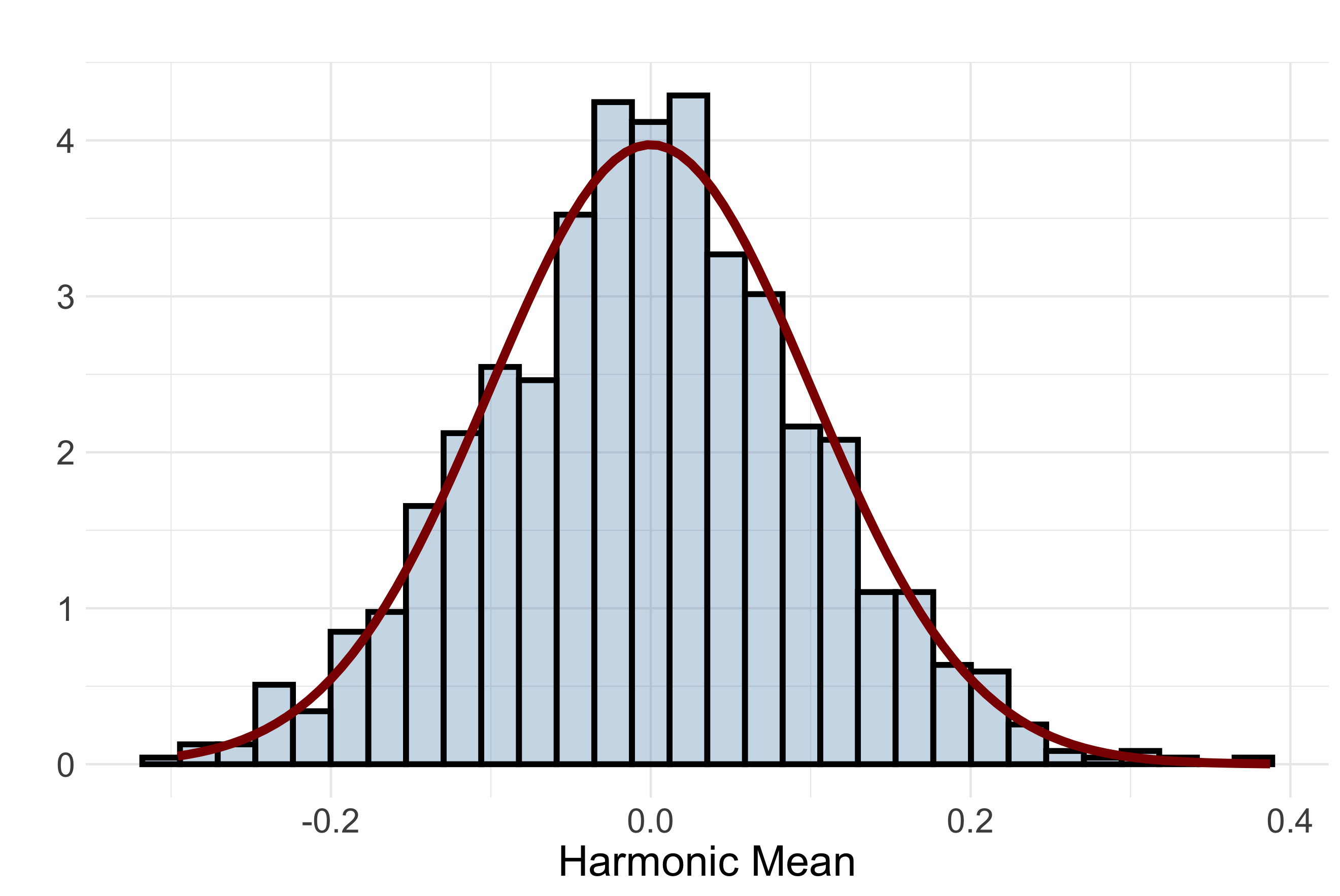}
  \caption{\label{sim1000}Monte Carlo simulation results for the universal central limit theorem with $n = 1\,000$.}
\end{figure}

\subsection{Numerical Experiments}
To illustrate the universal central limit theorem in practice we report the results of a simulation study where we sample 1\,000 datasets of size $n = 1\,000$, from the following models:
\begin{itemize}
\item Scenario A: LogNormal$(\mu, \sigma^2)$, with $(\mu, \sigma^2) = (2, 1)$. \vspace{-0.2cm}
\item Scenario B: Gamma$(a, b)$, with $(a, b) = (100, 1)$.\vspace{-0.2cm}
\item Scenario C: Uniform$(a, b)$, with $(a, b) = (1, 2)$.\vspace{-0.2cm}
\item Scenario D: Pareto$(\alpha)$, with $\alpha = 10$.
\end{itemize}
The results are reported in Figure~\ref{sim1000}. The left column corresponds to the particular case of the standard central limit theorem ($g(x) = x$), and which is included for comparison purposes. The middle and right columns illustrate the result for the particular cases of the geometric mean and harmonic mean (respectively $g(x) = \log x$ and $g(x) = 1/x$). As shown in the figure, the empirical distribution closely matches the theoretical prediction, aligning well with the asymptotic normal distribution implied by Theorem~\ref{UCLT}, with asymptotic variance given in~\eqref{avar}. 

\begin{figure}[H] \centering 
   \includegraphics[scale = .065]{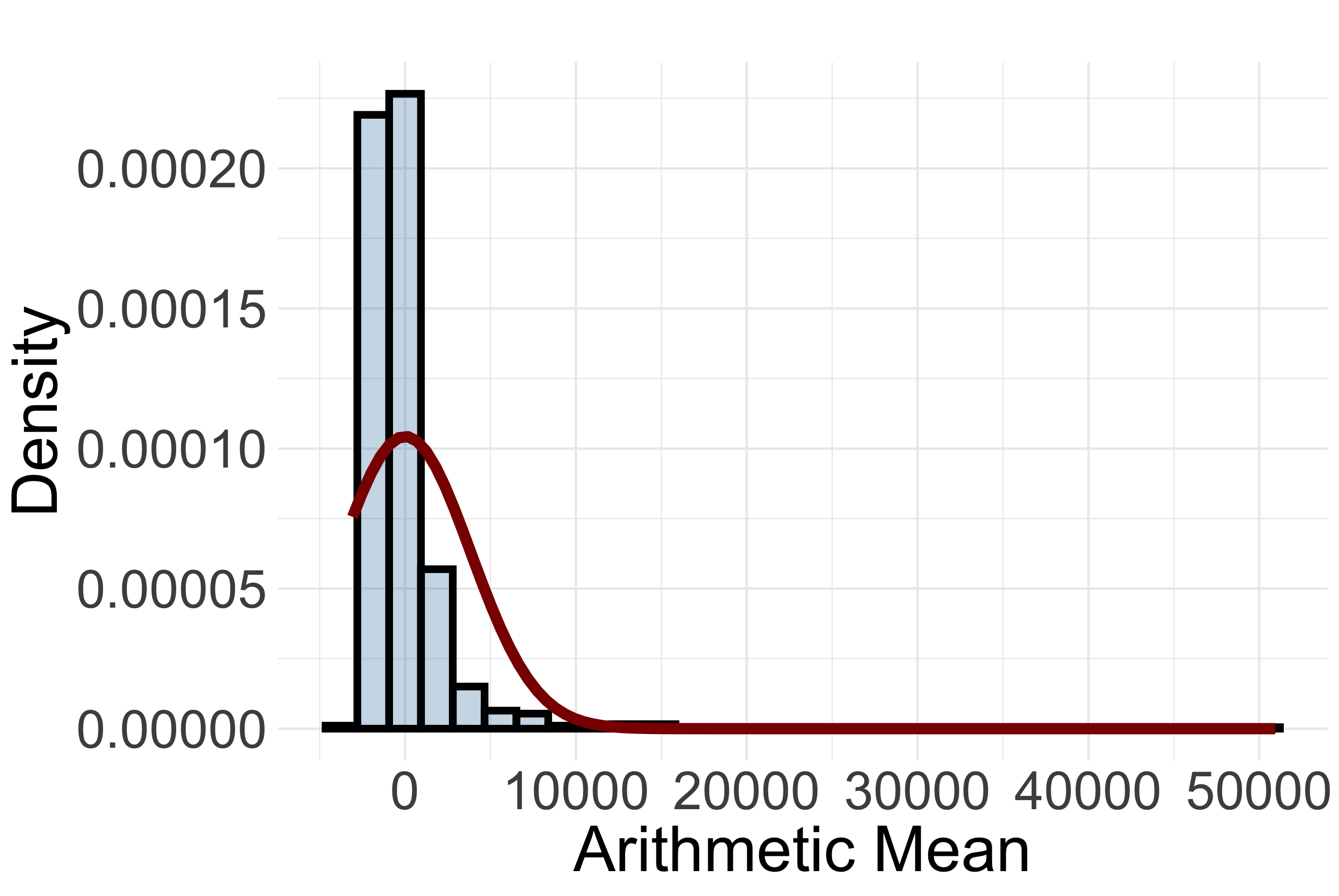}   \includegraphics[scale = .065]{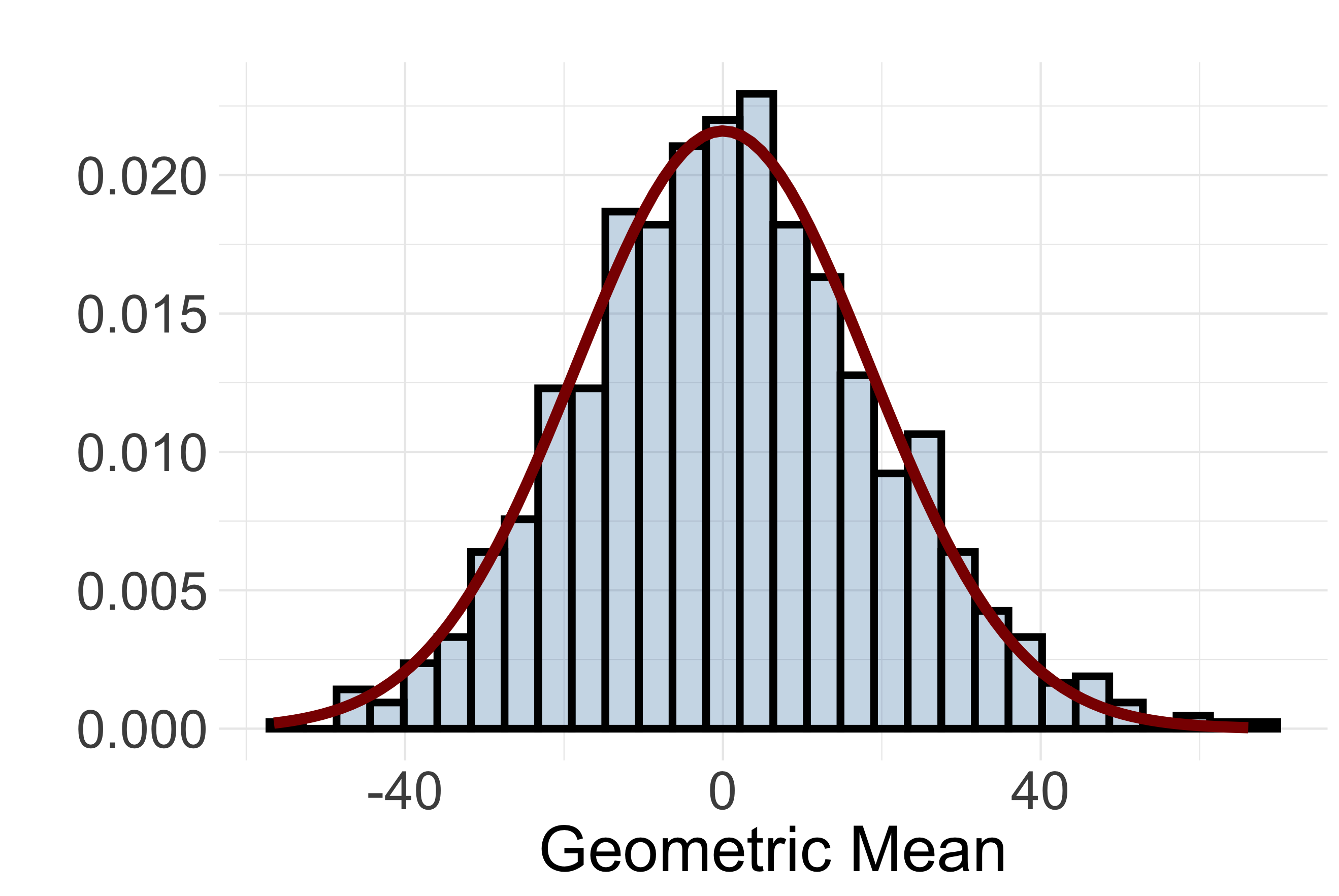} 
  \caption{\label{cltvs}Additional simulation results with $n = 1\,000$ from a $\text{LogNormal}(\mu, \sigma^2)$ distribution with $\mu = 2$ and $\sigma = 2.5$. As anticipated by the Edgeworth expansion in Section~\ref{ls}, convergence is slow for the arithmetic mean but fast for the geometric mean.}
\end{figure}

The theoretical insights from the Edgeworth expansion in Section~\ref{ls} are consistent with our numerical findings. For instance, in the LogNormal case the standard central limit theorem ($g(x)=x$) exhibits slow convergence when $\sigma$ is large, whereas convergence is considerably faster for the geometric mean ($g(x)=\log x$), since in that case $g(X)\sim \text{N}(\mu,\sigma^2)$ and hence, as predicted by the Edgeworth expansion, the normal approximation becomes accurate more rapidly. For an illustration of this fact see Figure~\ref{cltvs} which is based on LogNormal$(\mu, \sigma^2)$, with $\mu = 2$, and $\sigma = 2.5$. 

  \section{Proofs and Technical Details}
\subsection{Proof of Corollary~\ref{KLATB}}
Symmetry (i.e., A2) yields that 
  $f(x_1, \dots, x_n) = f(x_{\sigma_1}, \dots, x_{\sigma_n})$
for any permutation $\sigma$ of the indices $\{1, \dots, n\}$. This, together with Theorem~\ref{KLAT}, implies that for any permutation,
\begin{equation*}
  \sum_{j = 1}^{2 n + 1}l\left(\sum_{i = 1}^n \lambda_ig_j(x_i) \right) = \sum_{j = 1}^{2 n + 1}l\left(\sum_{i = 1}^n \lambda_ig_j(x_{\sigma_i}) \right).
\end{equation*}
In turn this implies that all $\lambda_i$ must be identical to ensure the equality above; that is, $\lambda_i = \lambda$ for all $i$. Hence, 
\begin{equation}\label{fin}
  \sum_{j = 1}^{2 n + 1}l\left(\lambda \sum_{i = 1}^n g_j(x_i) \right) =
  \sum_{j = 1}^{2 n + 1}l\left(\lambda \sum_{i = 1}^n g_j(x_{\sigma_i}) \right).
\end{equation}
Now, the only way for \eqref{fin} to hold for any permutation is if the $g_j$ are identical, as any variation would break symmetry; that is, $g_j = g$ for all $j$. This in turn implies that
\begin{equation*}
l\left(\lambda \sum_{i = 1}^n g(x_i) \right) =
l\left(\lambda \sum_{i = 1}^n g(x_{\sigma_i}) \right), 
\end{equation*}
Thus, we have proved that, under symmetry, both the inner and outer functions are constant, 
\begin{equation*}
  \begin{split}
    f(x_1, \dots, x_n) &= l\left(\lambda \sum_{i = 1}^n g(x_i) \right), 
  \end{split}
\end{equation*}
thereby establishing the final result. \qed
\subsection{Proof of Corollary~\ref{KLATC}}
First, note that by assumptions A1 and A2 $f$ is symmetric and continuous and hence Theorem~\ref{KLATC} can be applied. Next, it follows from Theorem~\ref{KLATC} that $g(x)$ is increasing and $\lambda > 0$, which implies that $$h(\mathbf{x}) \equiv \lambda \sum_{i = 1}^ng(x_i)$$ 
is increasing in each variable. Since by assumption A1 $f(\mathbf{x})$ is increasing in each variable, we next show that $y \mapsto l(y)$ has to be increasing, otherwise $f(\mathbf{x}) = l(h(\mathbf{x}))$ and $h(\mathbf{x})$ could not be both increasing in each variable. Indeed, note that monotonicity of $f$ and $h$ on each variable imply that for any $\varepsilon > 0$, 
$a \equiv h(x_1 + \varepsilon, \dots, x_n) > b \equiv h(x_1, \dots, x_n)$ 
and $l(a) > f(x_1 + \varepsilon, \dots, x_n) > f(x_1, \dots, x_n) = l(b)$.
Since the same argument applies for $x_2$, $x_3$, and so on, this proves that if $a > b$, then $l(a) > l(b)$,  hence establishing that $l$ is increasing and hence injective. Finally, combining A3 and A4 with \eqref{KLATB2} implies that for the case $n_0 = n$, $x = f(x_1, \dots, x_n)$, and thus since $l$ is injective,
\begin{equation*}
  \begin{split}
  & f(x_1, \dots, x_n) = f(x, \dots, x) \\ 
  &\Rightarrow l\left(\lambda \sum_{i = 1}^n g(x_i) \right) = l\left(\lambda \sum_{i = 1}^n g(f(x_1, \dots, x_n)) \right) \\ 
  &\Rightarrow \sum_{i = 1}^n g(x_i) = n \, g(f(x_1, \dots, x_n)) \\
  &\Rightarrow g^{-1}\left(\frac{1}{n} \sum_{i = 1}^n g(x_i)\right) = f(x_1, \dots, x_n),  
  \end{split}
\end{equation*}
and therefore the result is established. \qed
\subsection{Proof of Theorem~\ref{cont}}
Let $h \in C_{\text{LI}}(B)$ be arbitrary. Since $g$ is both strictly monotone and Lipschitz continuous, it can be easily shown that $g^{-1} \in C_{\text{LI}}(B)$ \citep[][Theorem~22.1]{estep2002}. A standard argument via the mean value theorem implies that 
\begin{equation}\label{MVT}
  \|g^{-1} - h^{-1}\|_{\infty} \leq \frac{1}{m}\|g - h\|_{\infty},
\end{equation}
where
$m = \min\{\inf_{x \in A}g'(x), \inf_{x \in A}h'(x)\} \neq 0$ as both
$g$ and $h$ are strictly monotone.

Let $L$ be the Lipschitz constant of $g^{-1}$, let $G_i = g(x_i)/n$ and $H_i = h(x_i)/n$, for all $i$. 
Equation~\eqref{MVT} along with the triangle inequality and the fact that $g^{-1} \in C_{\text{LI}}(B)$, yield that 
\begin{equation*}
  \begin{split}
  \|M_g - M_h\|_{\infty} &= \textstyle \|g^{-1}(\sum_{i = 1}^nG_i) - h^{-1}(\sum_{i = 1}^nH_i)\|_{\infty}\\
                         &\leq \textstyle \|g^{-1}(\sum_{i = 1}^nG_i) - g^{-1}(\sum_{i = 1}^nH_i)\|_{\infty} + \\
                         & \hspace{0.5cm}\textstyle \|g^{-1}(\sum_{i = 1}^nH_i) - h^{-1}(\sum_{i = 1}^nH_i)\|_{\infty} \\
                         &\leq \textstyle L \|\sum_{i = 1}^n(G_i - H_i)\|_{\infty} + \|g^{-1} - h^{-1}\|_{\infty} \\
                         &\leq \textstyle L \|g - h\|_{\infty} + 1/m\,\|g - h\|_{\infty} \\
                         &\leq \textstyle (L + 1 / m) \, \|g - h\|_{\infty}.
    \end{split}
  \end{equation*}
  Thus, to ensure that $\|M_g - M_h\|_{\infty} < \varepsilon$, for all $\varepsilon > 0$, it suffices to consider $\|g - h\|_{\infty} < \delta$, with $\delta = \varepsilon/(L + 1 / m)$. \qed

\section{Closing Remarks}
This article uncovers an elegant connection between two key results of Kolmogorov by showing that the class of regular means in~\eqref{rmean} can be derived directly from the celebrated superposition theorem; see Corollary~\ref{KLATC}. As a byproduct, this article highlights the often-overlooked power of the axiomatic perspective in Statistics and related fields. Assuming the continuity of a statistic $f(\mathbf{x})$ is equivalent to assuming that the superposition theorem holds. From this starting point, one can then explore which additional properties (e.g., symmetry) are desirable for a statistic $f(\mathbf{x})$ to possess, thereby further constraining its representation---bringing it closer to a concept of interest---in the spirit of Corollaries~\ref{KLATB} and~\ref{KLATC}.

A result on the smoothness of regular means is also established (Theorem~\ref{cont}), which can be viewed as a stability principle extending the continuity requirement imposed by A1. Finally, we provide insights into a universal central limit theorem, which may be regarded as a central limit theorem encompassing virtually all types of means, including the arithmetic, geometric, and harmonic means.

The version of the superposition theorem used above is sometimes also known as Kolmogorov--Arnol'd--Kahane--Lorentz--Sprecher theorem. Alternative non-monotone inner functions can be constructed \citep[][Theorem~4.2]{morris2021}. Indeed, another well-known version of the result removes the requirement of monotonicity in the inner functions while requiring $2n + 1$ outer functions. Specifically, let $f:A \to \mathbb{R}$ be a continuous function. Then it admits the representation
\begin{equation}\label{rep}
  f(\mathbf{x}) = \sum_{j=1}^{2n+1} \Phi_j^{(2)}\!\left(\sum_{i=1}^n \Phi_{j,i}^{(1)}(x_i)\right),
\end{equation}
for some continuous univariate functions $\Phi_{j,i}^{(1)}$ and $\Phi_{j}^{(2)}$, for all $i$ and $j$; see, for example, \citet{braun2009}. For our purposes, the version in Theorem~\ref{KLAT} is preferred, as it provides a more direct connection to regular means.

We close this paper with some notes on open problems and future directions. In my view, the full potential of the doubly additive representations in the Kolmogorov superposition theorem remains underappreciated in Statistics. As shown in this work, the theorem allows one to derive the entire class of regular means as well as other results (namely, Corollaries~\ref{KLATB} and~\ref{KLATC}).
More broadly, while not explored here, it is evident that \eqref{rep} has direct connections with statistical modeling. For instance, set $n = p$ and let $\mathbf{x} \in A \subset \mathbb{R}^p$ be a covariate vector. If all outer functions are proportional to the identity---that is, $\Phi_j^{(2)}(x) = x/{(2p+1)}$---and $\Phi^{(1)}_{j, i} = f_i$ the representation in~\eqref{rep} becomes additive, as in a generalized additive model (GAM) \citep{wood2006}
\begin{equation*}
  f(\mathbf{x}) = f_1^{(1)}(x_1) + \cdots + f_p^{(1)}(x_p).
\end{equation*}
From this perspective, models based in \eqref{rep} \citep[e.g.,][]{decarvalho2025, montanelli2020} can be viewed as doubly-additive neural models that naturally extend the class of GAMs. Finally, the potential of the superposition theorem in \eqref{rep} to lay the groundwork for flexible models for other objects of interest, such as joint densities and copulas, seems natural.

\bibliographystyle{asa2.bst}  
\bibliography{library.bib}       

\vspace{0.5cm} \footnotesize 
\noindent \textbf{Funding}: MdC is supported by the Leverhulme Trust, the Center for Investing Innovation, and Universidade de Aveiro---CIDMA, Funda\c c\~ao para a Ci\^encia e a Tecnologia (UID/4106/2025 and UID/PRR/4106/2025).

\end{document}